\documentclass[11pt]{article}

\usepackage{geometry}
\usepackage{amsmath}
\usepackage{amssymb}
\usepackage{amsthm}
\usepackage{graphicx}
\usepackage{subcaption}
\usepackage[round]{natbib}
\usepackage{hyperref}
\hypersetup{
	colorlinks=true,
	linkcolor=blue,
	citecolor=blue,
}
\usepackage{xcolor}
\usepackage{enumerate}
\usepackage{bm}
\usepackage{booktabs}
\usepackage{caption}

\theoremstyle{plain}

\newtheorem{proposition}{Proposition}
\newtheorem{lemma}{Lemma}

\title{
	A Note on ''Quasi-Maximum-Likelihood Estimation in Conditionally Heteroscedastic Time Series: A Stochastic Recurrence Equations Approach''
}
\author{
	\normalsize{Frederik Krabbe} \\
	\small{Aarhus University}
}
\date{}

\allowdisplaybreaks

\begin{document}
	
\maketitle

\begin{abstract}
	
	\cite{Bougerol1993} and \cite{StraumannMikosch2006} gave conditions under which there exists a unique stationary and ergodic solution to the stochastic difference equation $Y_t \overset{a.s.}{=} \Phi_t (Y_{t-1}), t \in \mathbb{Z}$ where $(\Phi_t)_{t \in \mathbb{Z}}$ is a sequence of stationary and ergodic random Lipschitz continuous functions from $(\textup{Y},|| \cdot ||)$ to $(\textup{Y},|| \cdot ||)$ where $(\textup{Y},|| \cdot ||)$ is a complete subspace of a real or complex separable Banach space. In the case where $(\textup{Y},|| \cdot ||)$ is a real or complex separable Banach space, \cite{StraumannMikosch2006} also gave conditions under which any solution to the stochastic difference equation $\hat{Y}_t \overset{a.s.}{=} \hat{\Phi}_t (\hat{Y}_{t-1}), t \in \mathbb{N}$ with $\hat{Y}_0$ given where $(\hat{\Phi}_t)_{t \in \mathbb{N}}$ is only a sequence of random Lipschitz continuous functions from $(\textup{Y},|| \cdot ||)$ to $(\textup{Y},|| \cdot ||)$ satisfies $\gamma^t || \hat{Y}_t - Y_t || \overset{a.s.}{\rightarrow} 0$ as $t \rightarrow \infty$ for some $\gamma > 1$. In this note, we give slightly different conditions under which this continues to hold in the case where $(\textup{Y},|| \cdot ||)$ is only a complete subspace of a real or complex separable Banach space by using close to identical arguments as \cite{StraumannMikosch2006}.
	
\end{abstract}

\section{Introduction}

Stochastic difference equations, such as
\begin{equation}
	Y_t \overset{a.s.}{=} \Phi_t (Y_{t-1}), \quad t \in \mathbb{Z}, \label{eq:SDE1}
\end{equation}
where $(\Phi_t)_{t \in \mathbb{Z}}$ is a sequence of stationary and ergodic random Lipschitz continuous functions from $(\textup{Y},|| \cdot ||)$ to $(\textup{Y},|| \cdot ||)$ where $(\textup{Y},|| \cdot ||)$ is a complete subspace of a real or complex separable Banach space $(\textup{B},|| \cdot ||)$, are ubiquitous in time series analysis. An example of a stochastic difference equation in time series analysis is the autoregressive (AR) model
\begin{equation*}
	Y_t \overset{a.s.}{=} \phi_0 + \phi_1 Y_{t-1} + \varepsilon_t, \quad t \in \mathbb{Z},
\end{equation*}
with $\phi_0 \in \mathbb{R}$ and $\phi_1 \in \mathbb{R}$ where $(\varepsilon_t)_{t \in \mathbb{Z}}$ is a sequence of independent and identically distributed (i.i.d.) real random variables with zero mean and variance $\sigma^2 > 0$. Another is the generalised autoregressive conditional heteroskedasticity (GARCH) model
\begin{equation*}
	\begin{aligned}
		Y_t &\overset{a.s.}{=} \sigma_t \varepsilon_t \\
		\sigma_t^2 &\overset{a.s.}{=} \omega + \alpha Y_{t-1}^2 + \beta \sigma_{t-1}^2
	\end{aligned}, \quad t \in \mathbb{Z},
\end{equation*}
with $\omega > 0$, $\alpha \geq 0$, and $\beta \geq 0$ where $(\varepsilon_t)_{t \in \mathbb{Z}}$ is a sequence of i.i.d. real random variables with zero mean and unit variance. Studying stochastic difference equations is therefore of great importance.

\cite{Bougerol1993} and \cite{StraumannMikosch2006} gave conditions under which there exists a unique stationary and ergodic solution to Equation \eqref{eq:SDE1} and any solution to the stochastic difference equation
\begin{equation}
	\bar{Y}_t \overset{a.s.}{=} \Phi_t (\bar{Y}_{t-1}), \quad t \in \mathbb{N}, \label{eq:SDE2}
\end{equation}
with $\bar{Y}_0$ given satisfies
\begin{equation}
	\left| \left| \bar{Y}_t - Y_t \right| \right| \overset{e.a.s.}{\rightarrow} 0 \quad \text{as} \quad t \rightarrow \infty; \label{eq:SDEC1}
\end{equation}
here, $\left| \left| X_t \right| \right| \overset{e.a.s.}{\rightarrow} 0$ as $t \rightarrow \infty$ means that $\gamma^{t} \left| \left| X_t \right| \right| \overset{a.s.}{\rightarrow} 0$ as $t \rightarrow \infty$ for some $\gamma > 1$.

Sometimes, however, the situation is more complicated than the one in Equation \eqref{eq:SDE2}. Consider, for instance, the following model
\begin{equation}
	Y_t \overset{a.s.}{=} \mu_t + \sigma_t \varepsilon_t, \quad t \in \mathbb{Z}, \label{eq:m1}
\end{equation}
where $(\varepsilon_t)_{t \in \mathbb{Z}}$ is a sequence of i.i.d. real random variables with zero mean and unit variance. Here,
\begin{align}
	\mu_t &\overset{a.s.}{=} \omega_{\mu} + \alpha_{\mu} Y_{t-1} + \beta_{\mu} \mu_{t-1}, \quad \quad \quad \quad \quad \ t \in \mathbb{Z}, \label{eq:m2} \\
	\sigma_t^2 &\overset{a.s.}{=} \omega_{\sigma} + \alpha_{\sigma} (Y_{t-1}-\mu_{t-1})^2 + \beta_{\sigma} \sigma_{t-1}^2, \quad t \in \mathbb{Z}, \label{eq:m3}
\end{align}
with $\omega_{\mu} \in \mathbb{R}$, $\alpha_{\mu} \in \mathbb{R}$, $\beta_{\mu} \in \mathbb{R}$, $\omega_{\sigma} > 0$, $\alpha_{\sigma} \geq 0$, and $\beta_{\sigma} \geq 0$. Consider now the initialised versions of Equations \eqref{eq:m2} and \eqref{eq:m3}
\begin{align}
	\bar{\mu}_t &\overset{a.s.}{=} \omega_{\mu} + \alpha_{\mu} Y_{t-1} + \beta_{\mu} \bar{\mu}_{t-1}, \quad \quad \quad \quad \quad \ t \in \mathbb{N}, \label{eq:m4} \\
	\bar{\sigma}_t^2 &\overset{a.s.}{=} \omega_{\sigma} + \alpha_{\sigma} (Y_{t-1}-\bar{\mu}_{t-1})^2 + \beta_{\sigma} \bar{\sigma}_{t-1}^2, \quad t \in \mathbb{N}, \label{eq:m5}
\end{align}
with $\bar{\mu}_0$ and $\bar{\sigma}_0^2$ given. For this model, the results by \cite{Bougerol1993} and \cite{StraumannMikosch2006} above can be used for Equations \eqref{eq:m2} and \eqref{eq:m4} since $\mu_t \overset{a.s.}{=} \Phi_t^{\mu} (\mu_{t-1}), t \in \mathbb{Z}$ and $\bar{\mu}_t \overset{a.s.}{=} \Phi_t^{\mu} (\bar{\mu}_{t-1}), t \in \mathbb{N}$ with $\bar{\mu}_0$ given where $\Phi_t^{\mu} (\mu) = \omega_{\mu} + \alpha_{\mu} Y_{t-1} + \beta_{\mu} \mu$. They can, however, not be used for Equations \eqref{eq:m3} and \eqref{eq:m5} since $\sigma_t^2 \overset{a.s.}{=} \Phi_t^{\sigma} (\sigma_{t-1}^2), t \in \mathbb{Z}$ where $\Phi_t^{\sigma} (\sigma^2) = \omega_{\sigma} + \alpha_{\sigma} (Y_{t-1}-\mu_{t-1})^2 + \beta_{\sigma} \sigma^2$ but $\bar{\sigma}_t^2 \overset{a.s.}{=} \bar{\Phi}_t^{\sigma} (\bar{\sigma}_{t-1}^2), t \in \mathbb{N}$ with $\bar{\sigma}_0^2$ given where $\bar{\Phi}_t^{\sigma} (\sigma^2) = \omega_{\sigma} + \alpha_{\sigma} (Y_{t-1}-\bar{\mu}_{t-1})^2 + \beta_{\sigma} \sigma^2$.

In the case where $(\textup{Y},|| \cdot ||)$ is a real or complex separable Banach space, \cite{StraumannMikosch2006} also gave conditions under which any solution to the stochastic difference equation
\begin{equation}
	\hat{Y}_t \overset{a.s.}{=} \hat{\Phi}_t (\hat{Y}_{t-1}), \quad t \in \mathbb{N}, \label{eq:SDE3}
\end{equation}
with $\hat{Y}_0$ given where $(\hat{\Phi}_t)_{t \in \mathbb{N}}$ is only a sequence of random Lipschitz continuous functions from $(\textup{Y},|| \cdot ||)$ to $(\textup{Y},|| \cdot ||)$ satisfies
\begin{equation}
	\left| \left| \hat{Y}_t - Y_t \right| \right| \overset{e.a.s.}{\rightarrow} 0 \quad \text{as} \quad t \rightarrow \infty. \label{eq:SDEC2}
\end{equation}
The condition that $(\textup{Y},|| \cdot ||)$ is a real or complex separable Banach space is, however, very restrictive; $([0,\infty),| \cdot |)$ is, for instance, not a real separable Banach space, so the result by \cite{StraumannMikosch2006} cannot be used for Equations \eqref{eq:m3} and \eqref{eq:m5}. In this note, we give slightly different conditions under which Equation \eqref{eq:SDEC2} continues to hold in the case where $(\textup{Y},|| \cdot ||)$ is only a complete subspace of a real or complex separable Banach space by using close to identical arguments as \cite{StraumannMikosch2006}. This result can, for instance, be used for Equations \eqref{eq:m3} and \eqref{eq:m5} since $([0,\infty),| \cdot |)$ is a complete subspace of the real separable Banach space $(\mathbb{R},| \cdot |)$.

The note is organised as follows. Section \ref{sec:results} collects the results, Section \ref{sec:proofs} collects the proofs, and the lemmata used in the proofs are collected in Section \ref{sec:lemmata}.

\section{Results} \label{sec:results}

The first proposition, which is a special case of Theorem 3.1 in \cite{Bougerol1993} and Theorem 2.8 in \cite{StraumannMikosch2006}, gives conditions under which there exists a unique stationary and ergodic solution to Equation \eqref{eq:SDE1}. In the following,
\begin{equation*}
	\Lambda (\Phi_t) = \sup_{\underset{x \neq y}{x,y \in \textup{Y}}} \frac{\left| \left| \Phi_t (x) - \Phi_t (y) \right| \right|}{\left| \left| x - y \right| \right|}.
\end{equation*}
\begin{proposition} \label{prop:1}
	Assume that
	\begin{enumerate}[(i)]
		\item there exists a $y \in \textup{Y}$ such that $\mathbb{E} \left[ \log^{+} \left| \left| \Phi_0 (y) - y \right| \right| \right] < \infty$ and
		\item $\mathbb{E} \left[ \log^{+} \Lambda (\Phi_0) \right] < \infty$ and there exists an $r \in \mathbb{N}$ such that
		\begin{equation*}
			-\infty \leq \mathbb{E} \left[ \log \Lambda (\Phi_0^{(r)}) \right] < 0,
		\end{equation*}
		where
		\begin{equation*}
			\Phi_t^{(r)} (\cdot) = \Phi_t \circ \Phi_{t-1} \circ \cdots \circ \Phi_{t-r+1} (\cdot).
		\end{equation*}
	\end{enumerate}
	Then,
	\begin{equation*}
		Y_t \overset{a.s.}{=} \lim_{n \rightarrow \infty} \Phi_t^{(n)} (y), \quad t \in \mathbb{Z},
	\end{equation*}
	is the only stationary and ergodic solution to Equation \eqref{eq:SDE1}.
\end{proposition}

The next one, which is also a special case of Theorem 3.1 in \cite{Bougerol1993}, shows that Equation \eqref{eq:SDEC1} also holds.
\begin{proposition} \label{prop:2}
	Assume that the assumptions in Proposition \ref{prop:1} hold. Then,
	\begin{equation*}
		\left| \left| \bar{Y}_t - Y_t \right| \right| \overset{e.a.s.}{\rightarrow} 0 \quad \text{as} \quad t \rightarrow \infty.
	\end{equation*}
\end{proposition}

The following proposition, which is a generalisation of Theorem 2.10 in \cite{StraumannMikosch2006}, gives conditions under which Equation \eqref{eq:SDEC2} holds.
\begin{proposition} \label{prop:3}
	Assume that the assumptions in Proposition \ref{prop:1} hold and that
	\begin{enumerate}[(i)]
		\item $\mathbb{E} \left[ \log^{+} \left| \left| Y_0 \right| \right| \right] < \infty$,
		\item $\left| \left| \hat{\Phi}_t (y) - \Phi_t (y) \right| \right| \overset{e.a.s.}{\rightarrow} 0$ as $t \rightarrow \infty$ for the $y \in \textup{Y}$ in Proposition \ref{prop:1}, and
		\item $\Lambda (\hat{\Phi}_t-\Phi_t) \overset{e.a.s.}{\rightarrow} 0$ as $t \rightarrow \infty$.
	\end{enumerate}
	Then,
	\begin{equation*}
		\left| \left| \hat{Y}_t - Y_t \right| \right| \overset{e.a.s.}{\rightarrow} 0 \quad \text{as} \quad t \rightarrow \infty.
	\end{equation*}
\end{proposition}

\section{Proofs} \label{sec:proofs}

\subsection{Proof of Proposition \ref{prop:1}}

Let $d$ be the metric induced by the norm $|| \cdot ||$. Then, $(\textup{B},d)$ is a complete separable metric space, so, since a subspace of a metric space is a metric space, $(\textup{Y},d)$ is a metric space, which is complete by assumption and separable since a subspace of a separable metric space is separable. The conclusion thus follows from Theorem 3.1 in \cite{Bougerol1993} (existence) and Theorem 2.8 in \cite{StraumannMikosch2006} (uniqueness).

% https://proofwiki.org/wiki/Subspace_of_Metric_Space_is_Metric_Space

% https://proofwiki.org/wiki/Subspace_of_Separable_Metric_Space_is_Separable

% https://proofwiki.org/wiki/Subspace_of_Complete_Metric_Space_is_Closed_iff_Complete

\subsection{Proof of Proposition \ref{prop:2}}

The conclusion follows from Theorem 3.1 in \cite{Bougerol1993}.

\subsection{Proof of Proposition \ref{prop:3}}

The conclusion follows if
\begin{equation}
	\left| \left| \hat{Y}_{tr} - Y_{tr} \right| \right| \overset{e.a.s.}{\rightarrow} 0 \quad \text{as} \quad t \rightarrow \infty. \label{eq:WTS}
\end{equation}
First, we show that
\begin{equation}
	\left| \left| \hat{\Phi}_t^{(n)} (y) - \Phi_t^{(n)} (y) \right| \right| \overset{e.a.s.}{\rightarrow} 0 \quad \text{as} \quad t \rightarrow \infty \label{eq:C2}
\end{equation}
and that
\begin{equation}
	\Lambda (\hat{\Phi}_t^{(n)}-\Phi_t^{(n)}) \overset{e.a.s.}{\rightarrow} 0 \quad \text{as} \quad t \rightarrow \infty \label{eq:C3}
\end{equation}
for all $n \in \mathbb{N}$ by induction. Equations \eqref{eq:C2} and \eqref{eq:C3} hold for $n = 1$ by assumption. Assume thus that Equations \eqref{eq:C2} and \eqref{eq:C3} hold for $n$. First,
\begin{align*}
	\left| \left| \hat{\Phi}_t^{(n+1)} (y) - \Phi_t^{(n+1)} (y) \right| \right| &= \left| \left| \hat{\Phi}_t \circ \hat{\Phi}_{t-1}^{(n)} (y) - \Phi_t \circ \Phi_{t-1}^{(n)} (y) \right| \right| \\
	& \leq \left| \left| \hat{\Phi}_t \circ \hat{\Phi}_{t-1}^{(n)} (y) - \hat{\Phi}_t \circ \Phi_{t-1}^{(n)} (y) \right| \right| \\
	& \quad + \left| \left| \left( \hat{\Phi}_t \circ \Phi_{t-1}^{(n)} (y) - \Phi_t \circ \Phi_{t-1}^{(n)} (y) \right) - \left( \hat{\Phi}_t (y) - \Phi_t (y) \right) \right| \right| \\
	& \quad + \left| \left| \hat{\Phi}_t (y) - \Phi_t (y) \right| \right| \\
	& \leq \Lambda (\hat{\Phi}_t) \left| \left| \hat{\Phi}_{t-1}^{(n)} (y) - \Phi_{t-1}^{(n)} (y) \right| \right| + \Lambda (\hat{\Phi}_t - \Phi_t) \left| \left| \Phi_{t-1}^{(n)} (y) - y \right| \right| \\
	& \quad + \left| \left| \hat{\Phi}_t (y) - \Phi_t (y) \right| \right| \\
	& \leq \Lambda (\hat{\Phi}_t - \Phi_t) \left| \left| \hat{\Phi}_{t-1}^{(n)} (y) - \Phi_{t-1}^{(n)} (y) \right| \right| + \Lambda (\Phi_t) \left| \left| \hat{\Phi}_{t-1}^{(n)} (y) - \Phi_{t-1}^{(n)} (y) \right| \right| \\
	&\quad + \Lambda (\hat{\Phi}_t - \Phi_t) \left| \left| \Phi_{t-1}^{(n)} (y) - y \right| \right| + \left| \left| \hat{\Phi}_t (y) - \Phi_t (y) \right| \right|.
\end{align*}Equation \eqref{eq:C2} thus also holds for $n+1$ by Lemmata \ref{lem:1} and \ref{lem:2} since $\mathbb{E} \left[ \log^{+} \left| \left| \Phi_0 (y) - y \right| \right| \right] < \infty$, $\mathbb{E} \left[ \log^{+} \Lambda (\Phi_0) \right] < \infty$, $\left| \left| \hat{\Phi}_t^{(m)} (y) - \Phi_t^{(m)} (y) \right| \right| \overset{e.a.s.}{\rightarrow} 0$ as $t \rightarrow \infty$ for all $m \in \{1,2,...,n\}$, and $\Lambda (\hat{\Phi}_t^{(m)}-\Phi_t^{(m)}) \overset{e.a.s.}{\rightarrow} 0$ as $t \rightarrow \infty$ for all $m \in \{1,2,...,n\}$. Moreover,
\begin{align*}
	\Lambda (\hat{\Phi}_t^{(n+1)}-\Phi_t^{(n+1)}) &= \Lambda (\hat{\Phi}_t \circ \hat{\Phi}_{t-1}^{(n)} - \Phi_t \circ \Phi_{t-1}^{(n)}) \\
	& \leq \Lambda (\hat{\Phi}_t \circ \hat{\Phi}_{t-1}^{(n)} - \hat{\Phi}_t \circ \Phi_{t-1}^{(n)}) + \Lambda (\hat{\Phi}_t \circ \Phi_{t-1}^{(n)} - \Phi_t \circ \Phi_{t-1}^{(n)}) \\
	& \leq \Lambda (\hat{\Phi}_t) \Lambda (\hat{\Phi}_{t-1}^{(n)} - \Phi_{t-1}^{(n)}) + \Lambda (\hat{\Phi}_t - \Phi_t) \Lambda (\Phi_{t-1}^{(n)}) \\
	& \leq \Lambda (\hat{\Phi}_t-\Phi_t) \Lambda (\hat{\Phi}_{t-1}^{(n)} - \Phi_{t-1}^{(n)}) + \Lambda (\Phi_t) \Lambda (\hat{\Phi}_{t-1}^{(n)} - \Phi_{t-1}^{(n)}) + \Lambda (\hat{\Phi}_t - \Phi_t) \Lambda (\Phi_{t-1}^{(n)})
\end{align*}
since $\Lambda$ is submultiplicative. Equation \eqref{eq:C3} thus also holds for $n+1$ by Lemmata \ref{lem:1} and \ref{lem:2} since $\mathbb{E} \left[ \log^{+} \Lambda (\Phi_0) \right] < \infty$ and $\Lambda (\hat{\Phi}_t^{(m)}-\Phi_t^{(m)}) \overset{e.a.s.}{\rightarrow} 0$ as $t \rightarrow \infty$ for all $m \in \{1,2,...,n\}$. We now prove that Equation \eqref{eq:WTS} holds. We have that
\begin{align*}
	\left| \left| \hat{Y}_{tr} - Y_{tr} \right| \right| &= \left| \left| \hat{\Phi}_{tr}^{(r)} (\hat{Y}_{(t-1)r}) - \Phi_{tr}^{(r)} (Y_{(t-1)r}) \right| \right| \\
	& \leq \left| \left| \hat{\Phi}_{tr}^{(r)} (\hat{Y}_{(t-1)r}) - \hat{\Phi}_{tr}^{(r)} (Y_{(t-1)r}) \right| \right| + \left| \left| \hat{\Phi}_{tr}^{(r)} (Y_{(t-1)r}) - \Phi_{tr}^{(r)} (Y_{(t-1)r}) \right| \right| \\
	& \leq \Lambda (\hat{\Phi}_{tr}^{(r)}) \left| \left| \hat{Y}_{(t-1)r} - Y_{(t-1)r} \right| \right| \\
	& \quad + \left| \left| \left( \hat{\Phi}_{tr}^{(r)} (Y_{(t-1)r}) - \Phi_{tr}^{(r)} (Y_{(t-1)r}) \right) - \left( \hat{\Phi}_{tr}^{(r)} (y) - \Phi_{tr}^{(r)} (y) \right)\right| \right| + \left| \left| \hat{\Phi}_{tr}^{(r)} (y) - \Phi_{tr}^{(r)} (y) \right| \right| \\
	& \leq \Lambda (\hat{\Phi}_{tr}^{(r)}) \left| \left| \hat{Y}_{(t-1)r} - Y_{(t-1)r} \right| \right| + \Lambda ( \hat{\Phi}_{tr}^{(r)} - \Phi_{tr}^{(r)} ) \left| \left| Y_{(t-1)r} - y \right| \right| + \left| \left| \hat{\Phi}_{tr}^{(r)} (y) - \Phi_{tr}^{(r)} (y) \right| \right| \\
	& \leq \Lambda (\hat{\Phi}_{tr}^{(r)}) \Lambda (\hat{\Phi}_{(t-1)r}^{(r)}) \left| \left| \hat{Y}_{(t-2)r} - Y_{(t-2)r} \right| \right| \\
	& \quad + \Lambda ( \hat{\Phi}_{tr}^{(r)} - \Phi_{tr}^{(r)} ) \left| \left| Y_{(t-1)r} - y \right| \right| + \left| \left| \hat{\Phi}_{tr}^{(r)} (y) - \Phi_{tr}^{(r)} (y) \right| \right| \\
	& \quad + \Lambda (\hat{\Phi}_{tr}^{(r)}) \left( \Lambda ( \hat{\Phi}_{(t-1)r}^{(r)} - \Phi_{(t-1)r}^{(r)} ) \left| \left| Y_{(t-2)r} - y \right| \right| + \left| \left| \hat{\Phi}_{(t-1)r}^{(r)} (y) - \Phi_{(t-1)r}^{(r)} (y) \right| \right| \right) \\
	& \ \, \vdots \\
	& \leq \left( \prod_{k=1}^{t} \Lambda (\hat{\Phi}_{kr}^{(r)}) \right) \left| \left| \hat{Y}_{0} - Y_{0} \right| \right| \\
	& \quad + \sum_{k=1}^{t} \left( \prod_{l=k+1}^{t} \Lambda (\hat{\Phi}_{lr}^{(r)}) \right) \left( \Lambda ( \hat{\Phi}_{kr}^{(r)} - \Phi_{kr}^{(r)} ) \left| \left| Y_{(k-1)r} - y \right| \right| + \left| \left| \hat{\Phi}_{kr}^{(r)} (y) - \Phi_{kr}^{(r)} (y) \right| \right| \right).
\end{align*}
First, we show that 
\begin{equation}
	\left( \prod_{k=1}^{t} \Lambda (\hat{\Phi}_{kr}^{(r)}) \right) \left| \left| \hat{Y}_{0} - Y_{0} \right| \right| \overset{e.a.s.}{\rightarrow} 0 \quad \textup{as} \quad t \rightarrow \infty. \label{eq:wts1}
\end{equation}
We have that
\begin{equation*}
	\left( \prod_{k=1}^{t} \Lambda (\hat{\Phi}_{kr}^{(r)}) \right) \left| \left| \hat{Y}_{0} - Y_{0} \right| \right| \leq \left( \prod_{k=1}^{t} \left( \Lambda (\Phi_{kr}^{(r)}) + \left| \Lambda (\hat{\Phi}_{kr}^{(r)}) - \Lambda (\Phi_{kr}^{(r)}) \right| \right) \right) \left| \left| \hat{Y}_{0} - Y_{0} \right| \right|,
\end{equation*}
where
\begin{equation*}
	\left| \Lambda (\hat{\Phi}_{kr}^{(r)}) - \Lambda (\Phi_{kr}^{(r)}) \right| \leq \Lambda (\hat{\Phi}_{kr}^{(r)}-\Phi_{kr}^{(r)}) \overset{e.a.s.}{\rightarrow} 0 \quad \text{as} \quad k \rightarrow \infty.
\end{equation*}
By the generalised monotone convergence theorem and Lemma \ref{lem:2},
\begin{equation*}
	\mathbb{E} \left[ \log \left(\Lambda (\Phi_{0}^{(r)}) + \varepsilon \right) \right] \downarrow \mathbb{E} \left[ \log \Lambda (\Phi_{0}^{(r)})\right] \quad \text{as} \quad \varepsilon \downarrow 0
\end{equation*}
since $\mathbb{E} \left[ \log^{+} \Lambda (\Phi_{0}) \right] < \infty$, so there exists an $\varepsilon_1 > 0$ such that
$\mathbb{E} \left[ \log \left(\Lambda (\Phi_{0}^{(r)}) + \varepsilon_1 \right) \right] < 0$ since $\mathbb{E} \left[ \log \Lambda (\Phi_{0}^{(r)})\right] < 0$. Thus, by Lemma \ref{lem:3},
\begin{equation*}
	\prod_{k=1}^{t} \left( \Lambda (\Phi_{kr}^{(r)}) + \left| \Lambda (\hat{\Phi}_{kr}^{(r)}) - \Lambda (\Phi_{kr}^{(r)}) \right| \right) \overset{e.a.s.}{\rightarrow} 0 \quad \text{as} \quad t \rightarrow \infty
\end{equation*}
since $\left| \Lambda (\hat{\Phi}_{kr}^{(r)}) - \Lambda (\Phi_{kr}^{(r)}) \right| \leq \varepsilon_1$ a.s. for all sufficiently large $k \in \mathbb{N}$ and $\mathbb{E} \left[ \log \left(\Lambda (\Phi_{0}^{(r)}) + \varepsilon_1 \right) \right] < 0$. Equation \eqref{eq:wts1} thus holds. We now show that also 
\begin{equation}
	\sum_{k=1}^{t} \left( \prod_{l=k+1}^{t} \Lambda (\hat{\Phi}_{lr}^{(r)}) \right) \left( \Lambda ( \hat{\Phi}_{kr}^{(r)} - \Phi_{kr}^{(r)} ) \left| \left| Y_{(k-1)r} - y \right| \right| + \left| \left| \hat{\Phi}_{kr}^{(r)} (y) - \Phi_{kr}^{(r)} (y) \right| \right| \right) \overset{e.a.s.}{\rightarrow} 0 \quad \text{as} \quad t \rightarrow \infty. \label{eq:wts2}
\end{equation}
First, by Lemmata \ref{lem:1} and \ref{lem:2}, 
\begin{equation*}
	\Lambda ( \hat{\Phi}_{kr}^{(r)} - \Phi_{kr}^{(r)} ) \left| \left| Y_{(k-1)r} - y \right| \right| + \left| \left| \hat{\Phi}_{kr}^{(r)} (y) - \Phi_{kr}^{(r)} (y) \right| \right| \overset{e.a.s.}{\rightarrow} 0 \quad \text{as} \quad k \rightarrow \infty
\end{equation*}
since $\mathbb{E} \left[ \log^{+} \left| \left| Y_{0} \right| \right| \right] < \infty$, $\left| \left| \hat{\Phi}_{kr}^{(r)} (y) - \Phi_{kr}^{(r)} (y) \right| \right| \overset{e.a.s.}{\rightarrow} 0$ as $k \rightarrow \infty$, and $\Lambda ( \hat{\Phi}_{kr}^{(r)} - \Phi_{kr}^{(r)} ) \overset{e.a.s.}{\rightarrow} 0$ as $k \rightarrow \infty$, so there exist a $0 < \gamma_1 < 1$ and a $C_1 \geq 0$ such that 
\begin{equation*}
	\Lambda ( \hat{\Phi}_{kr}^{(r)} - \Phi_{kr}^{(r)} ) \left| \left| Y_{(k-1)r} - y \right| \right| + \left| \left| \hat{\Phi}_{kr}^{(r)} (y) - \Phi_{kr}^{(r)} (y) \right| \right| \leq C_1 \gamma_1^{k} \quad a.s.
\end{equation*}
for all $k \in \mathbb{N}$. We thus have that
\begin{equation}
	\begin{aligned}
		& \sum_{k=1}^{t} \left( \prod_{l=k+1}^{t} \Lambda (\hat{\Phi}_{lr}^{(r)}) \right) \left( \Lambda ( \hat{\Phi}_{kr}^{(r)} - \Phi_{kr}^{(r)} ) \left| \left| Y_{(k-1)r} - y \right| \right| + \left| \left| \hat{\Phi}_{kr}^{(r)} (y) - \Phi_{kr}^{(r)} (y) \right| \right| \right) \\
		& \leq C_1 \sum_{k=1}^{t} \left( \prod_{l=k+1}^{t} \Lambda (\hat{\Phi}_{lr}^{(r)}) \right) \gamma_1^{k} \quad a.s.
	\end{aligned} \label{eq:proof1}
\end{equation}
Now, note that there exists an $\tilde{\varepsilon} > 0$ such that $\log \gamma_1 < \mathbb{E} \left[ \log \left( \Lambda (\Phi_{0}^{(r)}) + \tilde{\varepsilon} \right) \right] < 0$. Let $\tilde{\Lambda} (\hat{\Phi}_{lr}^{(r)}) = \Lambda (\hat{\Phi}_{lr}^{(r)}) + \tilde{\varepsilon}$ and $\tilde{\Lambda} (\Phi_{lr}^{(r)}) = \Lambda (\Phi_{lr}^{(r)}) + \tilde{\varepsilon}$. Then, continuing from Equation \eqref{eq:proof1}, 
\begin{equation}
	\begin{aligned}
		C_1 \sum_{k=1}^{t} \left( \prod_{l=k+1}^{t} \Lambda (\hat{\Phi}_{lr}^{(r)}) \right) \gamma_1^{k} &\leq C_1 \sum_{k=1}^{t} \left( \prod_{l=k+1}^{t} \tilde{\Lambda} (\hat{\Phi}_{lr}^{(r)}) \right) \gamma_1^{k} \\
		&= C_1 \sum_{k=1}^{t} \left( \prod_{l=1}^{t} \tilde{\Lambda} (\hat{\Phi}_{lr}^{(r)}) \right) \left( \prod_{l=1}^{k} \frac{\gamma_1}{\tilde{\Lambda} (\hat{\Phi}_{lr}^{(r)})} \right).
	\end{aligned} \label{eq:proof2}
\end{equation}
Note that
\begin{equation*}
	\prod_{l=1}^{k} \frac{\gamma_1}{\tilde{\Lambda} (\hat{\Phi}_{lr}^{(r)})} \leq \prod_{l=1}^{k} \left( \frac{\gamma_1}{\tilde{\Lambda} (\Phi_{lr}^{(r)})} + \left| \frac{\gamma_1}{\tilde{\Lambda} (\hat{\Phi}_{lr}^{(r)})} - \frac{\gamma_1}{\tilde{\Lambda} (\Phi_{lr}^{(r)})}\right| \right),
\end{equation*}
where, by the mean value theorem,
\begin{equation*}
	\left| \frac{\gamma_1}{\tilde{\Lambda} (\hat{\Phi}_{lr}^{(r)})} - \frac{\gamma_1}{\tilde{\Lambda} (\Phi_{lr}^{(r)})}\right| \leq \frac{\gamma_1}{\tilde{\varepsilon}^2} \left| \tilde{\Lambda} (\hat{\Phi}_{lr}^{(r)}) - \tilde{\Lambda} (\Phi_{lr}^{(r)}) \right| \overset{e.a.s.}{\rightarrow} 0 \quad \text{as} \quad l \rightarrow \infty.
\end{equation*}
By using the same arguments as above, there exists an $\varepsilon_2 > 0$ such that $\mathbb{E} \left[ \log \left( \frac{\gamma_1}{\tilde{\Lambda} (\Phi_{0}^{(r)})} + \varepsilon_2 \right) \right] < 0$ since $\mathbb{E} \left[ \log \frac{\gamma_1}{\tilde{\Lambda} (\Phi_{0}^{(r)})} \right] < 0$, so, by Lemma \ref{lem:3},
\begin{equation*}
	\prod_{l=1}^{k} \left( \frac{\gamma_1}{\tilde{\Lambda} (\Phi_{lr}^{(r)})} + \left| \frac{\gamma_1}{\tilde{\Lambda} (\hat{\Phi}_{lr}^{(r)})} - \frac{\gamma_1}{\tilde{\Lambda} (\Phi_{lr}^{(r)})}\right| \right) \overset{e.a.s.}{\rightarrow} 0 \quad \text{as} \quad k \rightarrow \infty
\end{equation*}
since $\left| \frac{\gamma_1}{\tilde{\Lambda} (\hat{\Phi}_{lr}^{(r)})} - \frac{\gamma_1}{\tilde{\Lambda} (\Phi_{lr}^{(r)})}\right| \leq \varepsilon_2$ a.s. for all sufficiently large $l \in \mathbb{N}$ and $\mathbb{E} \left[ \log \left( \frac{\gamma_1}{\tilde{\Lambda} (\Phi_{0}^{(r)})} + \varepsilon_2 \right) \right] < 0$. Thus, there exist a $0 < \gamma_2 < 1$ and a $C_2 \geq 0$ such that 
\begin{equation*}
	\prod_{l=1}^{k} \left( \frac{\gamma_1}{\tilde{\Lambda} (\Phi_{lr}^{(r)})} + \left| \frac{\gamma_1}{\tilde{\Lambda} (\hat{\Phi}_{lr}^{(r)})} - \frac{\gamma_1}{\tilde{\Lambda} (\Phi_{lr}^{(r)})}\right| \right) \leq C_2 \gamma_2^{k} \quad a.s.
\end{equation*}
for all $k \in \mathbb{N}$. Thus, continuing from Equation \eqref{eq:proof2},
\begin{equation*}
	C_1 \sum_{k=1}^{t} \left( \prod_{l=1}^{t} \tilde{\Lambda} (\hat{\Phi}_{lr}^{(r)}) \right) \left( \prod_{l=1}^{k} \frac{\gamma_1}{\tilde{\Lambda} (\hat{\Phi}_{lr}^{(r)})} \right) \leq \frac{C_1 C_2}{1-\gamma_2} \prod_{l=1}^{t} \tilde{\Lambda} (\hat{\Phi}_{lr}^{(r)}) \quad a.s.
\end{equation*}
Finally, we have that
\begin{equation*}
	\frac{C_1 C_2}{1-\gamma_2} \prod_{l=1}^{t} \tilde{\Lambda} (\hat{\Phi}_{lr}^{(r)}) \leq \frac{C_1 C_2}{1-\gamma_2} \prod_{l=1}^{t} \left( \tilde{\Lambda} (\Phi_{lr}^{(r)}) + \left| \Lambda (\hat{\Phi}_{lr}^{(r)}) - \Lambda (\Phi_{lr}^{(r)}) \right| \right),
\end{equation*}
where
\begin{equation*}
	\left| \Lambda (\hat{\Phi}_{lr}^{(r)}) - \Lambda (\Phi_{lr}^{(r)}) \right| \leq \Lambda (\hat{\Phi}_{lr}^{(r)}-\Phi_{lr}^{(r)}) \overset{e.a.s.}{\rightarrow} 0 \quad \text{as} \quad l \rightarrow \infty.
\end{equation*}
By using the same arguments as above, there exists an $\varepsilon_3 > 0$ such that $\mathbb{E} \left[ \log \left( \tilde{\Lambda} (\Phi_{0}^{(r)}) + \varepsilon_3 \right) \right] < 0$. Thus, by Lemma \ref{lem:3},
\begin{equation*}
	\prod_{l=1}^{t} \left( \tilde{\Lambda} (\Phi_{lr}^{(r)}) + \left| \Lambda (\hat{\Phi}_{lr}^{(r)}) - \Lambda (\Phi_{lr}^{(r)}) \right| \right) \overset{e.a.s.}{\rightarrow} 0 \quad \text{as} \quad t \rightarrow \infty
\end{equation*}
since $\left| \Lambda (\hat{\Phi}_{lr}^{(r)}) - \Lambda (\Phi_{lr}^{(r)}) \right| \leq \varepsilon_3$ a.s. for all sufficiently large $l \in \mathbb{N}$ and $\mathbb{E} \left[ \log \left( \tilde{\Lambda} (\Phi_{0}^{(r)}) + \varepsilon_3 \right) \right] < 0$. Equation \eqref{eq:wts2} thus also holds. This concludes the proof.

\section{Lemmata} \label{sec:lemmata}

\begin{lemma}[Lemma 2.1 in \cite{StraumannMikosch2006}] \label{lem:1}
	Let $(X_t)_{t \in \mathbb{N}_0}$ be a sequence of non-negative real random variables and let $(Y_t)_{t \in \mathbb{N}_0}$ be a sequence of stationary non-negative real random variables with $\mathbb{E} \left[ \log^{+} Y_0 \right] < \infty$. If $X_t \overset{e.a.s.}{\rightarrow} 0$ as $t \rightarrow \infty$, then
	\begin{equation*}
		X_t Y_t \overset{e.a.s.}{\rightarrow} 0 \quad \text{as} \quad t \rightarrow \infty.
	\end{equation*}
\end{lemma}

\begin{lemma}[Lemma 2.2 in \cite{StraumannMikosch2006}] \label{lem:2}
	Let $X$ and $Y$ be non-negative real random variables. Then,
	\begin{equation*}
		\mathbb{E} \left[ \log^{+} \left( X + Y \right) \right] \leq 2 \log 2 + \mathbb{E} \left[ \log^{+} X \right] + \mathbb{E} \left[ \log^{+} Y \right]
	\end{equation*}
	and
	\begin{equation*}
		\mathbb{E} \left[ \log^{+} \left( X Y \right) \right] \leq \mathbb{E} \left[ \log^{+} X \right] + \mathbb{E} \left[ \log^{+} Y \right].
	\end{equation*}
\end{lemma}

\begin{lemma}[Lemma 2.4 in \cite{StraumannMikosch2006}] \label{lem:3}
	Let $(X_t)_{t \in \mathbb{N}_0}$ be a sequence of stationary and ergodic non-negative real random variables with $- \infty \leq \mathbb{E} \left[ \log X_0 \right] < 0$. Then,
	\begin{equation*}
		\prod_{t=0}^{n} X_t \overset{e.a.s.}{\rightarrow} 0 \quad \text{as} \quad n \rightarrow \infty.
	\end{equation*}
\end{lemma}

\bibliography{references}

\end{document}